# The human brain and mathematics: recent research and educational questions

Gary E. Davis & Mercedes A. McGowen

University of Massachusetts Dartmouth, Dartmouth Massachusetts, USA
William Rainey Harper College, Palatine, Illinois, USA

Corresponding author: Gary E. Davis
Email: gdavis@umassd.edu
ORCID: 0000-0001-8668-4877

Second author: Mercedes A. McGowen
Email: mercmcgowen@sbcglobal.net
ORCID: 0000-0002-1070-359X

**Abstract**

*Recent experiments in neuroscience demonstrate the existence of a brain system that deals with mathematical thought and that is disjoint from the language areas of the brain. We provide an overview of these and related mathematically activated brain regions and raise a number of questions of interest to mathematics educators, including the nature of the brain's functioning in proof production, the role of compressible mental units in mathematical thought, and the nature of attention and working memory.*

## Introduction

New understandings of the functioning of human brains engaged in mathematics raise interesting questions for mathematics educators. Novel lines of research are suggested by neuroscientific findings, and new light is shed on some longstanding issues in mathematics education.

At the basis of a growing neuroscientific understanding of the mathematical functioning of human brains is the delineation of specific brain regions recruited and active during mathematical activity. Understanding these brain regions and their mathematical functionality brings a certain clarity and foundation for outstanding mathematics education concerns.

Mathematics educators have recently begun focusing much more on the nature of the human brain's mathematical functioning (although in the case of David Tall that focus goes back over 42 years). Examples are: De Smedt et al. 2010; De Smedt & Grabner 2015; De Smedt & Verschaffel 2010; Devlin 2010; Grabner et al. 2010; Kieran 2017; McGowen & Davis 2019; Peters & De Smedt 2018; Sellars 2018; Tall 1978, 1994, 2000, 2004, 2019; Verschaffel, Lehtinen & Van Dooren 2016). It seems timely, therefore, to help strengthen a growing connection between mathematics educators interested in, and neuroscientists carrying out empirical work on, the mathematical functioning of the human brain.

In this overview we focus mainly on the recent structural and functional findings of Stanislas Dehaene and colleagues (Amalric & Dehaene 2016, 2018, 2019; Dehaene et al. 2003; Dehaene 2009, 2010, 2011; Dehaene et al. 1999; Izard & Dehaene 2008; Maruyama et al. 2012; Pinheiro-Chagas et al. 2108; Sablé-Meyer et al. 2020). This is not to dismiss the empirical work of other neuroscientists actively researching aspects of the human brain's mathematical functioning (for example: Butterworth 1999; Hubbard et al. 2005; Kadosh & Walsh 2009; Kaufmann, Kucian & von Aster 2015; Krueger et al. 2008; Lee at al. 2007; Menon 2010, 2015; Rickard et al. 2000). Rather, we find the body of work of Dehaene and colleagues sufficiently broad and rich in detail to provide mathematics educators with a reliable up to date road map of what is known about different brain structures and functions relevant to mathematical thinking. That road map, we suggest, is helpful in formulating testable hypotheses relevant to mathematics education.

We discuss below issues that to us are important theoretically or practically for mathematics educators, issues that to our minds have particular relevance for mathematics education, and what the neuroscientific findings seem to indicate about such matters.

**Mathematics and language**

Connections between mathematical development and language use have been a topic of considerable theorizing and empirical research in the mathematics education community: see, for example, Bickmore-Brand 1990; Chapin et al. 2009; Clarke, Waywood & Stephens 1993; DeJarnette & González 2016; Ellerton & Clarkson 1996; Hunting 1998; Morgan 2006; O'Halloran 2015; Moschkovich 2015; Powell et al. 2017; Schleppegrell 2010; Shield & Galbraith 1998; and Wilkinson 2018, 2019.

When we are engaged in learning or discussing mathematics we read and talk with others, and in communicating mathematical thought we talk, write and draw diagrams. We might expect therefore, as many mathematics education researchers have either postulated or assumed, there is a tight connection between the human brain's mathematical functioning and its utilization of language.

The recent findings of Amalric & Dehaene (2016, 2018, 2019) have highly significant implications for the widely assumed connection between mathematics and language. In a nutshell, Amalric & Dehaene find there are brain regions, distinctly different from the language areas of the brain, that are recruited by expert mathematicians when reading mathematically-based text. Moreover the brain's language areas, in the experiments they describe, were largely *inactive* during processing of mathematical text by expert mathematicians.

This suggests, at least for expert mathematicians, a strong separation in the brain for mathematics on the one hand, and language on the other: quite the opposite scenario to that often postulated in the mathematics education literature. What are the findings of Amalric & Dehaene, and what are the mathematically-functioning brain regions activated in their experiments? Their major findings are summarized below.

In one experiment (Amalric & Dehaene 2016), 15 professional mathematicians and 15 humanities specialists (all French researchers or college level teachers) were presented with spoken mathematical and non-mathematical statements. These statements all had high level content: mathematical statements covered algebra, analysis, topology and geometry, while non-mathematical statements dealt with knowledge of nature and history.

Participants were given 4 seconds to judge if individual statements were true or false. Their brains were scanned during this procedure using functional MRI imaging. Amalric & Dehaene found a clear separation of brain areas used to interpret the mathematical statements between mathematicians and non-mathematicians:

> "… we show that high-level mathematical reasoning rests on a set of brain areas that do not overlap with the classical left-hemisphere regions involved in language processing or verbal semantics. Instead, all domains of mathematics we tested (algebra, analysis, geometry, and topology) recruit a bi-lateral network, of prefrontal, parietal, and inferior temporal regions, which is also activated when mathematicians or non-mathematicians recognize and manipulate numbers mentally. Our results suggest that high-level mathematical thinking makes minimal use of language areas and instead recruits circuits initially involved in space and number." (Amalric & Dehaene 2016, p. 4909)

In a further study, Amalric & Dehaene (2019) utilized a group of 14 professional mathematicians who were given 2.5 seconds to decide if orally presented mathematical or non-mathematical statements (some of which are shown below) were true or false.
They found even simple mathematical statements and rote algebraic expressions did not invoke language areas of the brains of participants, while invoking specific non-language areas as observed in their earlier experiments. Examples of mathematical and non-mathematical statements used by Amalric & Dehane (2019) are:

**Mathematical**

- The sine function is periodical
- $(a+b)(a–b) = a^2 – b^2$
- $(x-1)(x+1) = x^2 – 1$
- $\sin(x+3\pi/2) = -\cos(x)$
- The section of a sphere by a plane is always a point
- Lp spaces are separable
- Some matrices are diagonalizable
- Hyperboloids are not connected
- Some order relations are not transitive
- Some infinite sets are not countable

**Non-mathematical**

- London buses are red
- Rock'n roll is a musical style characterize by a slow tempo
- The Paris metro was built before the Istanbul one
- Some ocean currents are warm
- Orange blossoms are not perfumed
- Some plants are not climbing
- Tigers are not carnivores

What are the brain regions activated during processing of mathematical statements by professional mathematicians, to which Amalric & Dehaene (2016, 2018, 2019) refer? Even for simple mathematical statements there were two general areas of the brain activated in their interpretation by professional mathematicians: intraparietal sulci in parietal lobe and inferior temporal gyrus on temporal lobe. With more complex mathematical statements, activation occurred in dorsal prefrontal cortex, a functional area of prefrontal cortex, and when explicitly geometric statements were processed by expert mathematicians, activation occurred in calcarine sulci, where primary visual cortex is located.

(1) *Intraparietal sulci*, on parietal lobe (one for each hemisphere):

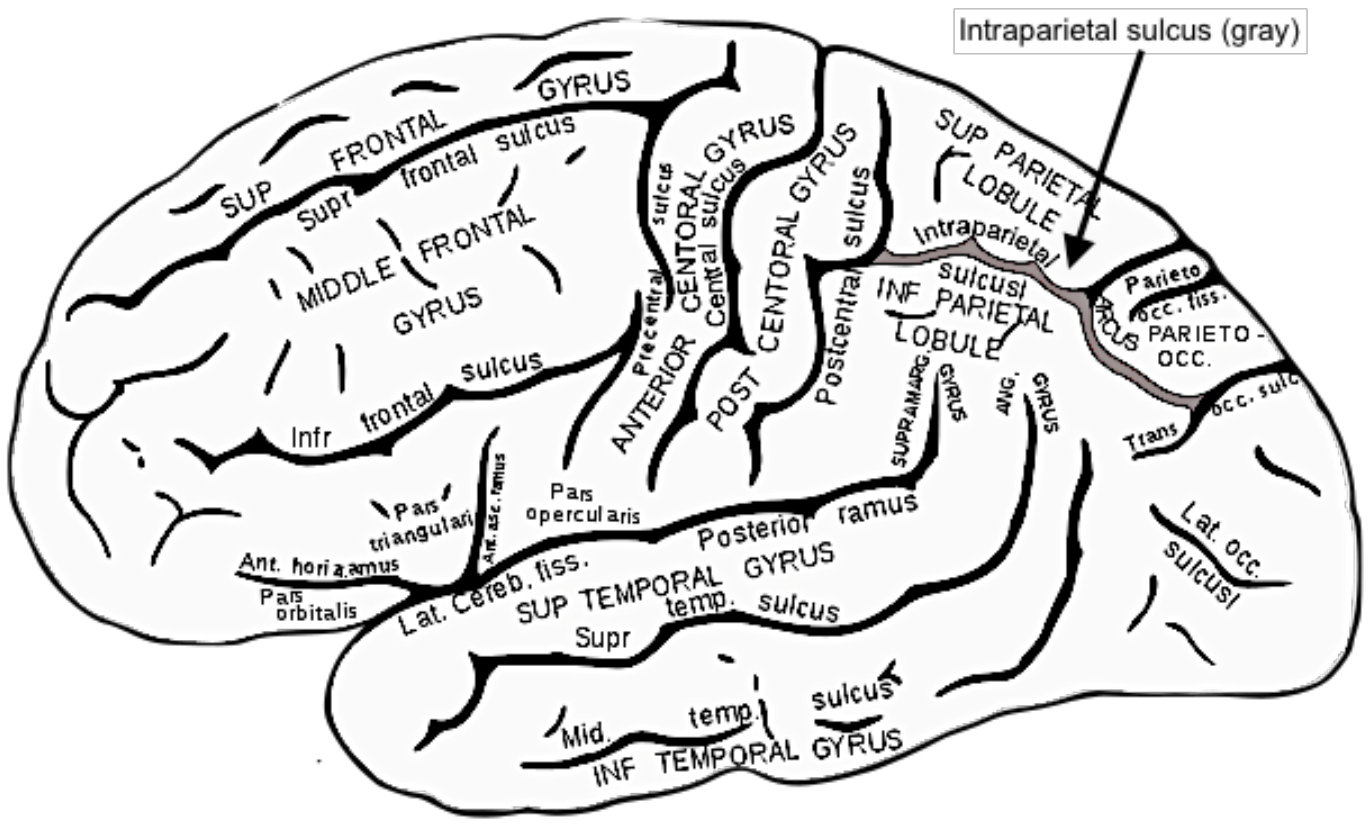

**Fig 1** Intraparietal sulcus on parietal lobe (top right: gray. Figure modified from work in public domain at https://commons.wikimedia.org/wiki/File:Gray726.svg

Apart from the activation of intraparietal sulci during interpretation of even simple mathematical statements by professional mathematicians, as observed by Amalric & Dehaene, intraparietal sulci, which contain several functionally distinct subregions (ventral, lateral, medial, anterior), are implicated in perceptual-motor coordination and visual attention – for example, directing eye movements, pointing, reaching, grasping, and object manipulation (Culham & Kanwisher 2001). Intraparietal sulci are also implicated in processing of symbolic numerical information (Canton, et al., 2006), and in visuospatial working memory (Todd & Marois 2004). Dormal & Pesenti (2009) found, in a study of fourteen right-handed French-speaking males (mean age: 21 ± 2.3 years), engaged in estimating numerosity and spatial magnitude, that left intraparietal sulcus was involved in processing numerosity and not spatial magnitude, while right intraparietal sulcus was engaged in representation of both spatial and numerical magnitude. The situation might be different for left-handed people, yet this study shows convincingly the roles of intraparietal sulci in processing of numerical and spatial magnitude. Hawes et al. (2019) studied brain regions active in mental rotation, arithmetic, and symbolic number processing and found that intraparietal sulci were "the largest and most consistent region of overlap across all three cognitive tasks" (p. 39), with left intraparietal sulcus being the largest region of activation for symbolic number and arithmetic, and right intraparietal sulcus the largest region of activation for mental rotation.

From Amalric & Dehaene's recent experiments it is not clear if a particular *subregion* (ventral, lateral, medial, anterior) of intraparietal sulcus is activated by professional mathematicians' processing of mathematical statements. Much of the research on functioning of intraparietal sulci has been carried out on monkey brains, however Culham & Kanwisher (2001) provide a summary of what is known about human intraparietal sulci and parietal lobes more generally. Gottlieb (2007) suggests that *lateral* interparietal sulci, which are known to be heavily involved in spatial and visual attention (Cohen 2014, p. 306) "may mediate interactions between spatial orienting and higher-level, more abstract, cognitive functions" (p. 14). Colby, Duhamel & Goldberg (1993) suggest that *ventral* subregions of intraparietal areas, which contain a predominance of neurons responsive to direction and speed, functions in analysis of visual motion. Colby & Goldberg (1999) concur as to the responsiveness of neurons in ventral intraparietal areas to vision and pointed out functioning of these areas in touch of the head and mouth. In the same article Colby & Goldberg point out the role of *medial* intraparietal areas in responding to a variety of stimuli within reaching distance of an individual, and that neurons in *anterior* intraparietal area are responsive to visual stimuli that an individual (monkey) is able to manipulate. In relation to *lateral* intraparietal areas, Colby & Goldberg summarize neural activity in these subregions of parietal lobes as contributing to spatial attention and spatial memory, in which only significant signals stimulate neuronal activity.

Generally, as Colby & Goldberg (1999) point out, intraparietal areas, and parietal lobes more generally, are involved in multiple representations of space, and in visual attention. Hindsight is a wonderful thing and looking back to the Colby & Goldberg article one might have imagined from their exposition that parietal lobes, and intraparietal sulci, particularly, might play a significant role in the interpretation of mathematical statements. That this is the case, and that there is such a striking dissociation form the language areas of the brain, is a truly significant finding of Amalric & Dehaene.

(2) *Inferior temporal gyrus*, on the temporal lobe:

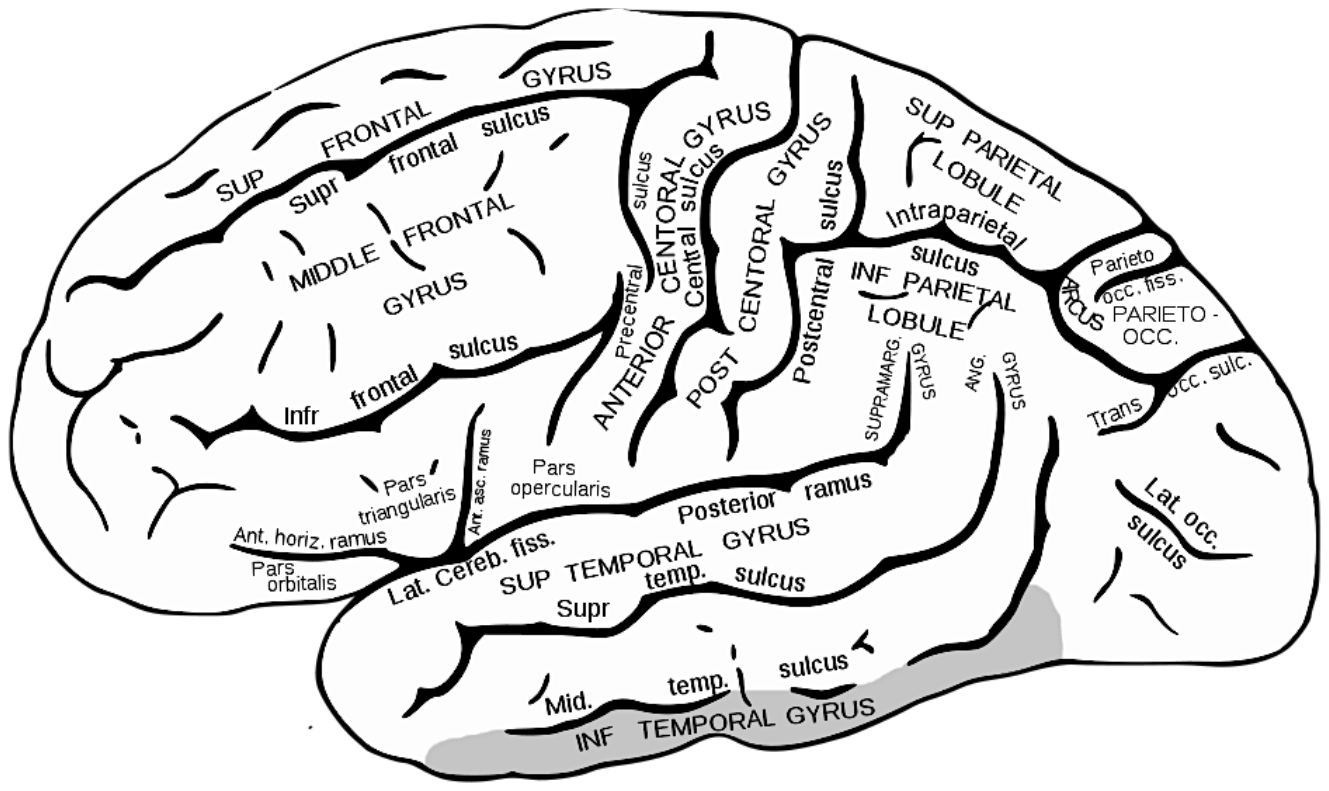

**Fig 2** Inferior temporal gyrus in temporal lobe (bottom: gray). Figure modified from work in public domain at https://commons.wikimedia.org/wiki/File:Gray726.svg

This region of temporal lobe is involved in visual processing, particularly the representation of objects, places, faces, and colors (Lafer-Sousa & Conway 2013), and in the recognition of numerals (Shum et al., 2013). Grotheer, Jeska & Grill-Spector (2018) and Pinheiro-Chagas, et al. (2108) also detail the role of posterior inferior temporal gyrus in numerical processing. Miyashita (1993) identifies inferior temporal lobe as a region where visual perception connects with memory and imagery. DiCarlo, Zoccolan & Rust (2012) provide evidence that the ability to rapidly recognize and categorize objects resides ultimately, after feedforward information from other brain regions, in neuronal representation in inferior temporal cortex.

Intraparietal sulci and inferior temporal gyri were activated in *all* Amalric & Dehaene's experiments when professional mathematicians were presented with mathematical statements, simple or complex, to interpret. With more complex mathematical statements Amalric & Dehaene (2016) also found activation of dorsal prefrontal cortex:

(3) *Dorsal prefrontal cortex* is a functional (not anatomical) area of prefrontal cortex and lies in the lateral part of Brodmann's areas 9 and 46.

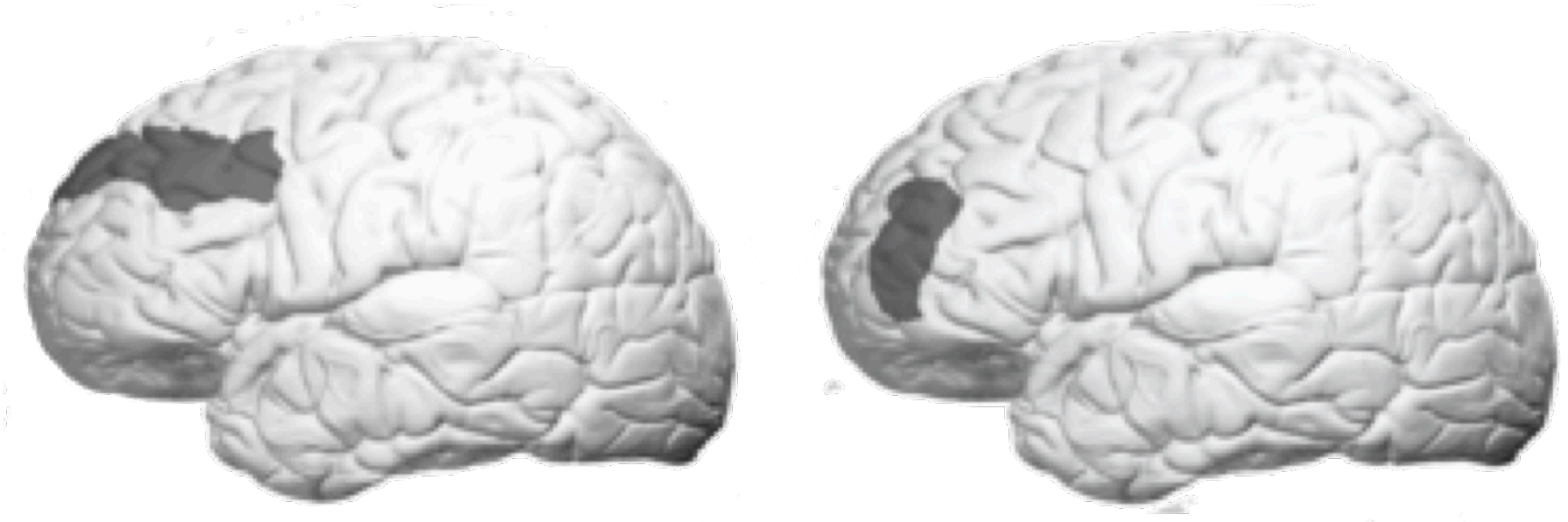

**Fig 3** Brodmann areas 9 (left) and 46 (right). Images modified from works under a Creative Commons CC BY-SA 2.1 JP license at https://commons.wikimedia.org/wiki/File:Brodmann_area_9_lateral.jpg and https://commons.wikimedia.org/wiki/File:Brodmann_area_46_lateral.jpg

Dorsolateral prefrontal cortex is involved, but not exclusively, in the management of working memory and cognitive flexibility (Kaplan, Gimbel & Harris 2016) as well as planning, and abstract reasoning (Bettcher 2017). Prefrontal cortex undergoes massive changes during adolescence, including pruning of neuronal connections and a growth in systems that increase dopamine which regulates the intensity of goal-directed behaviors (Spear 2000; Puglisi-Allegra & Ventura 2012). Adolescence is generally taken to include the years 12 – 18, and potentially up to age 25 (Spear, 2000). Given Amalric & Dehaene's findings that expert mathematicians' dorsolateral prefrontal cortex is activated by interpretation of more complex mathematical statements, and given the massive ongoing changes in prefrontal cortex in adolescence, it seems to us to be a matter of considerable importance for the mathematical education of adolescents to ascertain more accurately just what the changes in adolescent prefrontal cortex entail for student interpretation of complex mathematical statements.

Additionally, when more overtly geometric statements, such as "The angle between i and 1+i equals $\pi/4$", or "Any equilateral triangle can be divided into two right triangles", and trigonometric statements, were presented to the mathematicians studied by Amalric & Dehaene, further non-language brain regions were activated. Specifically, questions about complex numbers and trigonometry activated calcarine sulcus, where primary visual cortex is located. The posterior part of calcarine sulcus activates for the central visual field and the anterior part for the peripheral visual field:

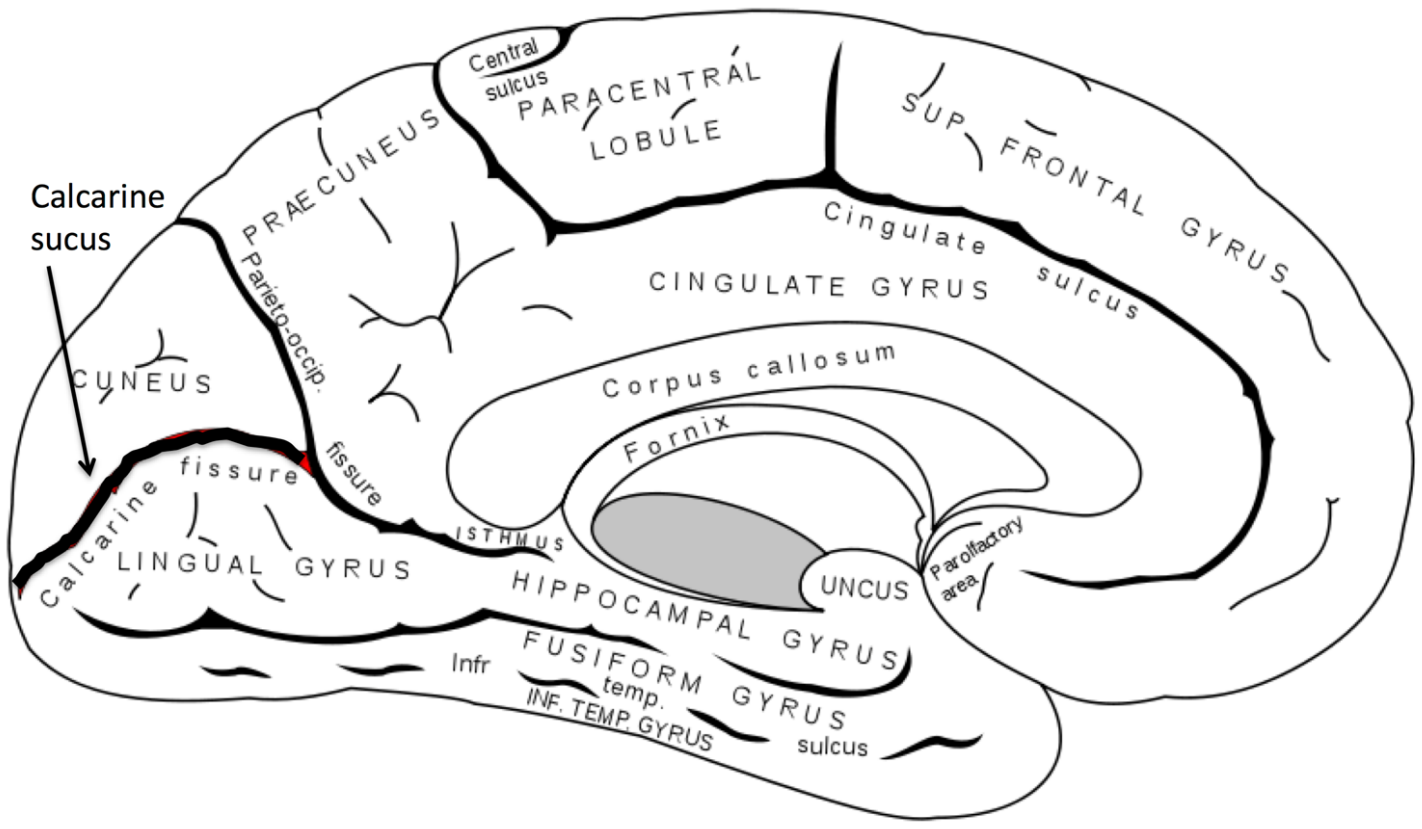

**Fig 4** Calcarine sulcus. Figure modified from a work in the public domain at https://commons.wikimedia.org/wiki/File:Gray727_calcarine_sulcus.svg

*Numerical and algebraic concepts and language*

Further indicators of neural dissociation of processing of mathematical syntax and classical language areas of the brain are reported in several other studies:

- Gelman & Butterwoth (2005) argue, on the basis of theory and data, that numerical concepts have neural bases that are independent of language.
- Monti, Parsons & Osherson (2012) found a neural dissociation between algebra and natural language.
- Maruyama et al. (2008) demonstrated that language areas were not recruited when students were asked to process nested algebraic expressions.

*Developmental issues*

Although the work of Amalric & Dehaene we report here deals with adult, expert, mathematicians, similar dissociations between mathematics and language have been observed, in several studies, in infants and young children: see Almaric & Dehaene (2018, section 3(f): Child development).

## The brain and proof

*Proof comprehension*

The experiments of Amalric & Dehaene all deal with interpretation of written or spoken mathematical language. This leads us to believe that when mathematicians are reading a mathematical proof, or listening to an explanation of a mathematical argument, specific non-language regions of their brains will be activated, especially intraparietal sulci and inferior temporal gyri. However, just as interpretation of more geometric mathematical statements, especially those couched in terms of complex numbers and trigonometry, led to activation of calcarine sulcus in the primary visual cortex, it is possible that the peculiar nature of interpreting a mathematical proof may involve other non-language areas of the brain, or, indeed, may invoke activation of language areas of the brain:

- Which brain regions are activated when professional mathematicians, and mathematics students, are engaged in attempting to understand a written or spoken proof. Are language areas of the brain more heavily involved in proof understanding?

We note that fMRI may not be the most appropriate imaging method to address these issues since proof comprehension is generally extended over a lengthy time period, and fMRI is generally unsuitable for temporal resolution.

*Proof production*

We turn now to what we see as one of the great mysteries raised by Amalric & Dehaene's (2016, 2018, 2019) findings. Students, at various levels of their education, are exposed to the need to present, orally, or in writing, a mathematical argument as to why some observation is true. As simple examples, encountered even in elementary school, consider the following observations:

- A long rectangle is formed by chaining 1×1 squares:

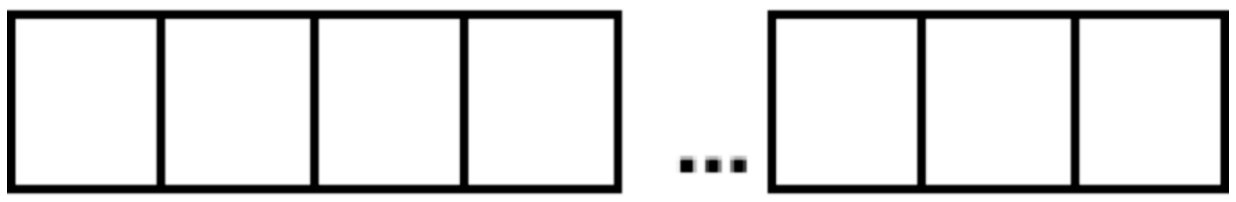

**Fig 5** A rectangle made from abutting unit squares

A table of the perimeter of the rectangle as a function of the number of squares is constructed:

| Number of squares | Perimeter of rectangle |
|---|---|
| 1 | 4 |
| 2 | 6 |
| 3 | 8 |
| 4 | 10 |
| 5 | 12 |
| 6 | 14 |

Table 1. Record of the perimeter of the rectangle as a function of the number of squares

A student says: “The perimeter is going up by 2 each time.” Why is that?
Another student says: “The perimeter is two more than twice the number of squares.” Why is that?

- Students are building polyominos from unit squares, of various shapes, some convex, some not. Some students observe the perimeter of every polyomino seen is even and conjecture this is always true. Why is that all polyominos have even perimeter?
- When we build all block towers from two colors of blocks, of a given height the number of block towers is always a power of 2:

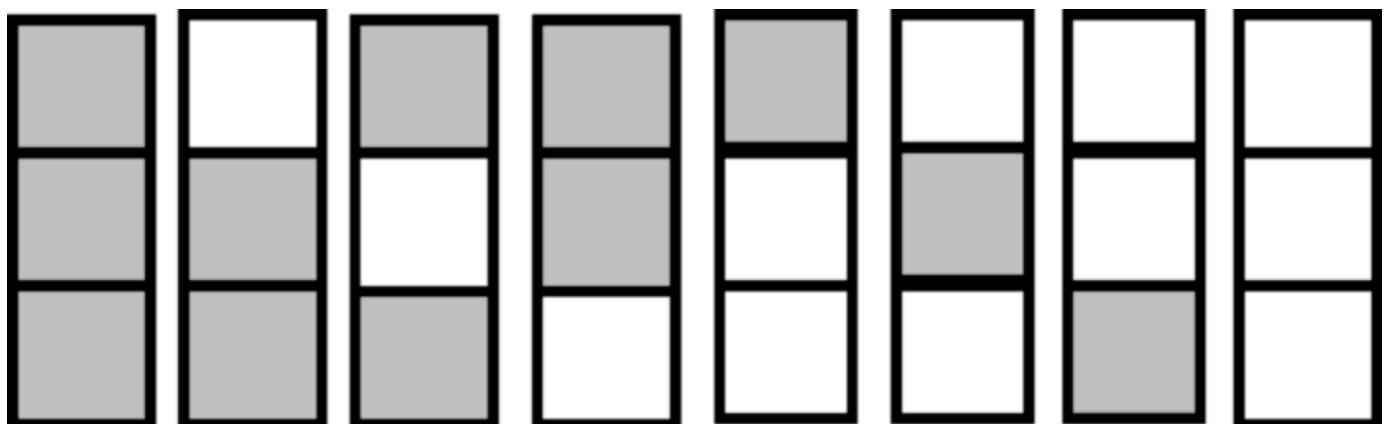

$2^3 = 8$ block towers from 2 colors, of height 3

**Fig 6** Block towers, of height 3 built using one or two colors

Why is that?

- In the game of “Make 21” two players take turns in writing down either the number 1 or the number 2. At each step a cumulative total of the numbers written is kept. The first player to reach 21 wins. Does player 1 or player 2 have a winning strategy and why?

How might we ascertain what brain region and functions are activated as students grapple with providing spoken and written explanations for their reasons why? fMRI imaging has poor temporal resolution, seemingly necessary for the extended times that may be required to produce explanations. However, fMRI might tell us something about which brain regions and functions are activated in these circumstances for expert mathematicians: that is, if we provide simple, elementary problems such as those above, previously unseen, to expert mathematicians, they are very likely to produce oral explanations quite rapidly. Those explanations, by their very nature, involve spoken language, so we expect to see the language areas of expert mathematician’s brains activated. For a more complex, yet still elementary, proof we might pose the question of why, when the edges of a complete graph on 6 vertices are colored either red or blue there must be a monochrome triangle in the graph. Assuming an expert mathematician, or a student, has not already seen this problem the reason is likely to be worked out over many minutes, requiring something other than fMRI imaging to determine where and when brain regions and functions are activated. Another complex, yet still elementary, question is why, in a 5×5 grid it is not possible to walk from the indicated square – marked with an x – to every other square once and only once by walking up, down, left or right to an adjacent square:

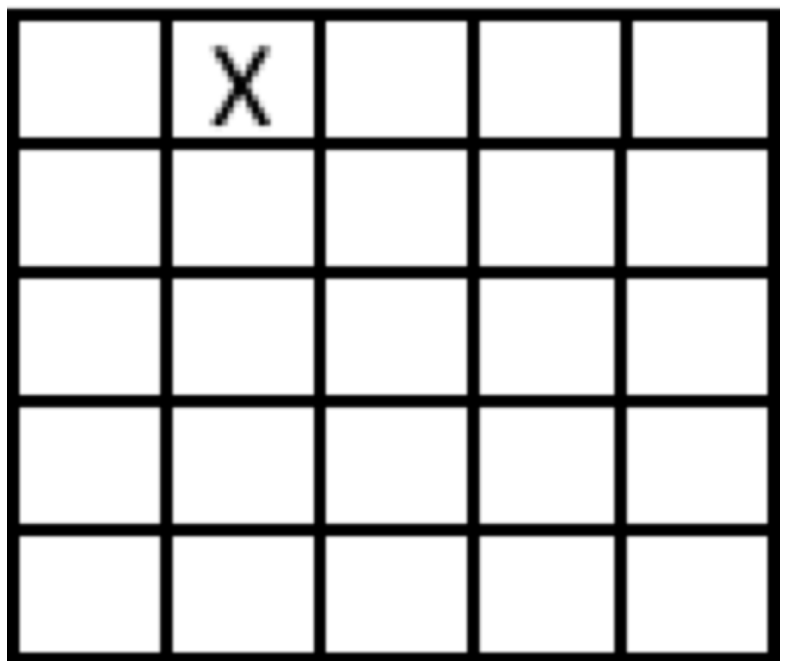

**Fig 7** Possible to walk to horizontally or vertically adjacent squares from the square marked, and cover each square exactly once?

Proof, by its very nature, involves expressing mathematical thought in words and mathematical signs, so one would imagine that it necessarily involves activity in the language areas of the human brain. Suppose, for example, that a student in a real analysis course is asked to write a proof that the sequence of rational numbers $p(n)$ defined recursively as:

$$p(1) = 1 \text{ and } p(n+1) = 1 - \frac{p(n)}{2}$$

for $n \geq 1$, converges to $\frac{2}{3}$.

A proof in the standard undergraduate ϵ-N style might go as follows:

Let $\varepsilon > 0$ and choose an integer $N > 1/(3\varepsilon)$. For each $n \geq 1$ we have

$$\left|p(n+1) - \frac{2}{3}\right| = \left|1 - \frac{p(n)}{2} - \frac{2}{3}\right| = \left|\frac{1}{3} - \frac{p(n)}{2}\right|$$

$$= ½\left|\frac{2}{3} - p(n)\right| = ½\left|p(n) - \frac{2}{3}\right|$$

so, by induction, $\left|p(n) - \frac{2}{3}\right| = \frac{1}{3\times 2^{n-1}}$ .Then, for $n > N$ we have

$$2^{n-1} > 2^{N-1} = (1+1)^{N-1} \geq 1 + (N-1) = N > 1/(3\varepsilon)$$

so $\left|p(n) - \frac{2}{3}\right| = \frac{1}{3\times 2^{n-1}} < \varepsilon$ .

On the face of it this sort of written proof involves elements of ordinary language, a lot of mathematical signs, an apparently tight linguistic structure, and a forward plan. Typically, students find constructing this sort of proof hard (Weber 2001). For us, a question that needs to be teased apart is: do students find it hard because their mathematical ideas are occurring not in the language areas of their brains but because, and perhaps by necessity, they must use the language areas of their brains to express those thoughts in writing? In other words, is a large problem – but far from the only problem – in student proof-writing due to a translation difficulty from areas of the brain dealing with mathematical thought to the language areas of the brain? While we may be far from being able to address this question yet, it is central to understanding the difficulty mathematics students have with more formal proof.

There remains the intriguing possibility that during proof production there is a continual to-and-fro feedback between the regions of the brain Amalric and Dehaene identify as being active in mathematical comprehension, on the one hand, and the brain's language centers on the other. Most mathematicians, and some students, will have experienced the phenomenon of putting thoughts down in writing, involving words and mathematical signs, and then looking at those words and signs to stimulate further thought – the next step in a proof. Thompson (2013) writes:

> … it can be unsettling to realize how much thinking already happens outside our skulls. … The physicist Richard Feynman once got into an argument about this with the historian Charles Weiner. Feynman understood the extended mind; he knew that writing his equations and ideas on paper was crucial to his thought. But when Weiner looked over a pile of Feynman's notebooks, he called them a "wonderful record of his day-to-day work." No, no, Feynman replied testily. They weren't a record of his thinking process. They *were* his thinking process:
> "I actually did the work on the paper," he said.
> "Well," Weiner said, "the work was done in your head, but the record of it is still here."
> "No, it's not a *record*, not really. It's *working*. You have to work on paper and this is the paper. Okay?" (p. 7)

Our classroom experience is that there is indeed a to-and-fro feedback between mathematical thought and linguistic expression of that thought, in spoken or written language. One of us recalls an incident in an undergraduate real analysis class in which a student, Stephanie, was listening to an explanation of a fellow student, both spoken and written on a white board. Stephanie was waving her left hand and apparently pointing at something out of field, while trying to speak to the presenting student. The instructor interrupted and asked Stephanie if she had something to say. She replied she did, but could not find the words, while still waving her left hand. The instructor said it seemed she was seeing something in her mind as she waved her hand. Stephanie was adamant that she could see what she wanted to communicate, she just could not find the words. Our experience is that this phenomenon is common, and we believe it highly likely that there remains to be discovered the neural mechanisms by which activity in the above-described regions of the brain are translated to the language areas of the brain for oral or written communication, and how that communication then stimulates further mathematical thought. This complex set of interactions is what happens in mathematics classrooms, and it appears that neuroscience may be on the verge of telling us in greater detail what is happening in the brains of students during these interactions. That clearer, deeper, knowledge of brain function hopefully will help reflective mathematics educators to better assist students.

**Neural basis of compression**

William Thurston remarked on the human brain's capacity for compressing mathematical thought:

> "Mathematics is amazingly compressible: you may struggle a long time, step by step, to work through some process or idea from several approaches. But once you really understand it and have the mental perspective to see it as a whole, there is often a tremendous mental compression. You can file it away, recall it quickly and completely when

> you need it, and use it as just one step in some other mental process. The insight that goes with this compression is one of the real joys of mathematics." (Thurston 1990, p. 848)

We have no reason to doubt Thurston, especially since our own experiences of learning mathematics are in accord with his statement. The question for us here is: what brain structures and functioning correspond to this "tremendous mental compression"?

Tall, in particular, has utilized Thurston's idea of mental compression and linked it to a feature of human brains engaged in mathematical thinking (Gray & Tall 2007; Tall 2007 2109b, 2020a, 2020b). Gray & Tall (2007, p. 23) write:

> "This paper considers mathematical abstraction arising through a natural mechanism of the biological brain in which complicated phenomena are compressed into thinkable concepts. The neurons in the brain continually fire in parallel and the brain copes with the saturation of information by the simple expedient of suppressing irrelevant data and focusing only on a few important aspects at any given time."

This aspect of compressibility in the brain, described by Gray & Tall, is essentially statistical redundancy as discussed, for example, by Barlow (2001), which is thought to play a major role in the functioning of vision (Olshausen & Field 2000). Does the compression of which Thurston writes constitute a different form of compression, utilizing specific and identifiable brain regions and functions, or is it simply compression due to redundancy reduction when brains are engaged in mathematical thought? Or, as we feel is more likely, does it have little to do with redundancy reduction, but involves the brain in the formation of "mental units", a particular mathematical phenomenon that relates to Thurston's notion of compression, but not in an obvious way, to statistical redundancy? We feel that Thurston's compression of process and ideas is more akin to what we think of as the formation of a repeatable mental unit, which we elaborate below, and which, in the spirit of this guide, we hope might be linked to specific brain regions and functioning.

Let us give two examples. Consider the problem, posed to a class of year 10 students in the UK: "how many 19s are there in 97, and with what remainder?" The majority of the class attempted this problem via long division, realized they would already need to know the answer to begin, and began successively adding lots of 19. A very few students simply wrote "5, remainder 2" and when asked how they got that answer replied, in summary, that they saw 19 as almost 20, they knew 5 lots of 20 is 100, so coming back 5 gives 95: 2 less than 97. Essentially, they were, from our perspective (and not necessarily, or even, from theirs), able form a mental unit:

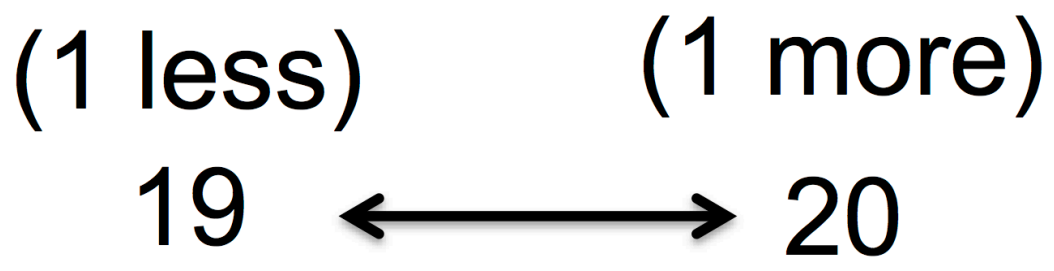

**Fig 8** A hypothesized mental unit

which they could, in thought, repeat 5 times

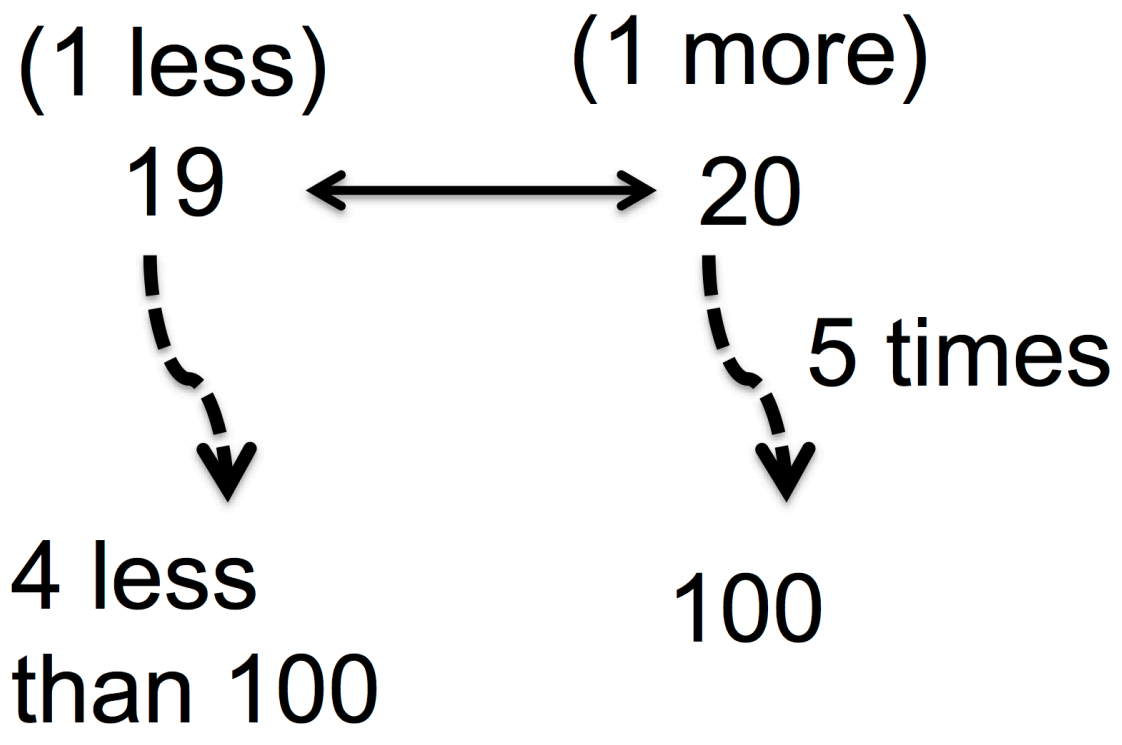

**Fig 9** Repetition of a mental unit

In the repetition, the connection between 19 and 20 – the essence of the unit – is preserved.

As a second example of a mental unit consider the problem of counting how many straight line segments there are in any given step of the construction of the Koch snowflake curve.

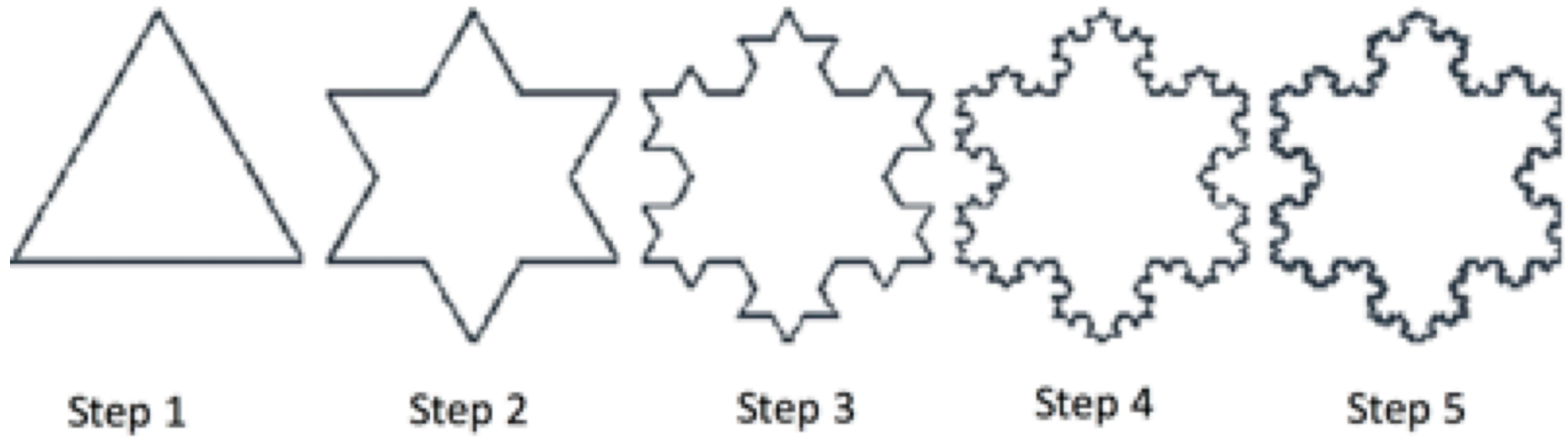

**Fig 10** Steps in the construction of the Koch snowflake curve

One student, in a group of undergraduates in an Australian university, was explaining how to obtain the number of straight line segments at step 4 of the snowflake construction by drawing the right half of the curve at that step and counting 1, 2, 3, 4, … . Another student commented they did not count, and when asked what they did replied, in summary, that at each stage of the construction, every line segment is replaced by 4 line segments as follows:

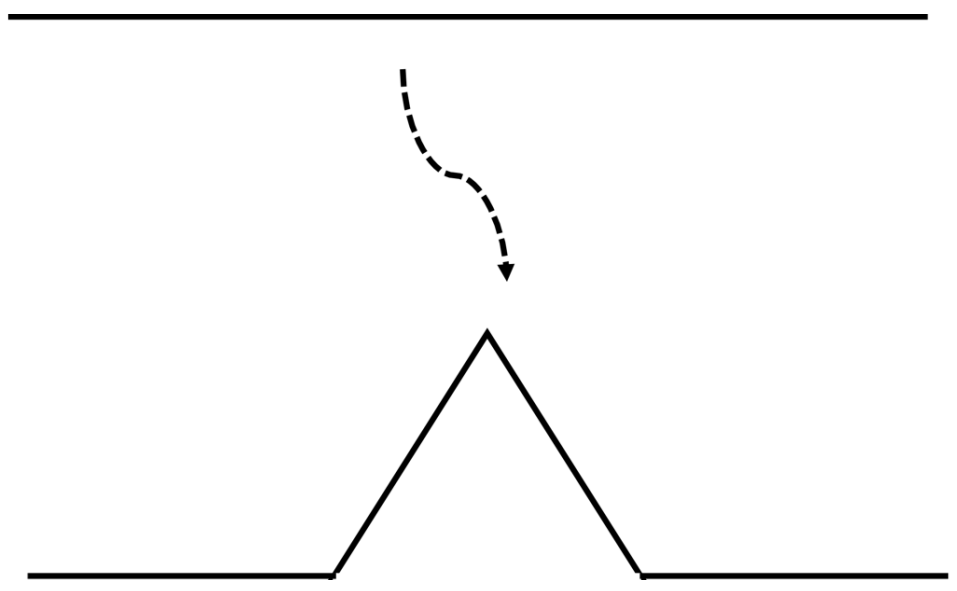

**Fig 11** A mental unit for the construction of the snowflake curve

Since there are 3 line segments at step 1, there are 3×4 = 12 at step 2, 12×4 = 48 at step 3, and 48×4 = 192 at step 4. Recognition of this repeatable mental unit allows the second student to bypass the laborious task of counting via one-one correspondence.

The formation of these – repeatable – mental units seems to us to be a form of compression that is distinctly different from redundancy reduction. Since this form of mental compression is so valuable in advancing mathematical thinking for individual students we ask: what is the neural basis of the formation of these repeatable mental units? In particular:

- What is happening in the brains of students at different ages and stages of development as they engage with problems such as these?
- What differences are there in recruitment of different brain regions and brain functioning between students who do utilize mental units like these and those who do not?
- Is there a common, identifiable collection of brain regions and functions involved in the formation of such mathematical units?
- Does the formation of mental units utilize the brain's visual processing regions, in particular, lateral intraparietal sulci and inferior temporal lobes?
- Do students at different developmental levels utilize essentially the same, or different, brain regions and functioning as and when they demonstrably form and use repeatable mental units, especially in basic arithmetical or geometric thinking?

**Focus of attention and working memory**

Students' initial focus of attention when working on a mathematical problem sets the scene for how they bring known acts and conditioned ways of working and thinking to bear on a problem. We discus focus of attention and working memory together, because in our pedagogic experience they are intimately linked and because they are also closely related neuroanatomically (LaBar et al. 1999).

Three brain regions are commonly activated in working memory: (a) parietal cortex, (b) dorsolateral prefrontal cortex and (c) supplementary motor cortex (Cohen 2014, p. 839). We note that the first two of these were found by Amalric & Dehaene to be activated during expert mathematician's processing of certain mathematical statements. This raises an obvious question: is a general working memory deficit associated with relatively poor mathematical performance *because* brain regions that (expert) mathematicians utilize, in interpreting mathematical statements, have functional deficits?

The neuropsychology of attention is a complex and detailed field (Cohen 2014), and there are specific contributions to attentional processes from the intraparietal sulci: brain regions activated in all Amalric & Dehaene's studies of interpretation of mathematical statements. Bisley & Goldberg (2006) suggest the lateral intraparietal area provides a *salience map*, whose function is to rank the significance of different parts of the visual field. The issue of salience is directly relevant to the work of two students described in McGowen & Davis (2019). In the context of a problem of modeling the dynamics of a toy rocket, one student focused on a quadratic model and quadratic regression to obtain parameters from given data, while another student focused on numerical values provided by data and attempted to use the quadratic formula for roots in an unproductive way. We suggest what a student finds salient in visual attention impacts what they bring to mind in working memory to address the problem to hand. How we might assist students to build stronger and more relevant salience maps is as yet unknown. However, the role of the intraparietal areas seems to be central to so doing.

There are at two forms of attention relevant to working memory: *executive attention* which is concerned with overall carrying out, including planning, of a real-time task, and *selective attention* which is invoked when attention during a task is oriented to specific aspects, objects or representations involved in the task to the exclusion of others. Vandierendonck (2014) argues that in many circumstances selective attention can produce an overload on working memory that inhibits executive attention. In other words, focus on a specific aspect of a task can interfere, in real time, with executive control of the execution of the task.

What attentional (and non-attentional) mechanisms contribute to keeping available information in working memory? This issue has been reviewed for verbal memory by Camos & Barrouillet (2014). Their conclusion is that just two systems sustain *verbal* information in working memory:

> "…we suggest that only two systems sustain the maintenance of information at the short term, counteracting the deleterious effect of temporal decay and interference.. A non-attentional mechanism of verbal rehearsal, similar to the one described by Baddeley in the phonological loop model, uses language processes to reactivate phonological memory traces. Besides this domain-specific mechanism, an executive loop allows the reconstruction of memory traces through an attention-based mechanism of refreshing." (p. 900)

The corresponding knowledge for mathematical information seems to be entirely unknown, and, given the major findings of Amalric & Dehaene on a separate brain system for processing mathematics, it seems likely that there might be distinct neurological aspects to the maintenance of mathematical information in working memory.

**Conclusions and discussion**

Recent brain imaging studies appear to be converging on two strongly related aspects of the human brain engaged in mathematical thought:

1. Processing of mathematical statements are, in the main, carried out by regions of the brain that are distinct from the language areas of the brain, and the latter are largely silent during such activity.
2. Several of the regions and functional areas of the brain that are activated during processing of mathematical statements are strongly involved in visual attention, visual processing, recall and working memory.

As a ball-park simplification, therefore, one might say that mathematics is not, in the main, language-based, but rather utilizes the brain's capacity for visual attention and visual representation. Since mathematical thought now appears to utilize heavily the human brain's visual attention and visual processing regions, we might speculate that the origins of mathematical thought are essentially visual in nature, and that, from an educational perspective, focusing on visual representations might therefore link productively with the brain's mathematical functioning. Amalric, Denghien & Deheane (2018, p. 314) studied 3 blind mathematicians who utilized the same brain mechanisms as previously studied mathematicians, in the absence of visual experience. They point out that "…the mechanisms by which formal mathematics emerges from proto-mathematical systems for numbers and space remain unknown. A possibility is that mathematical representations are rooted in visuospatial thinking and develop through visual experience."

In particular, a focus on visual representations might help us begin to understand the role and functioning of a student's initial focus of attention and links to recall and memory. We suggest a student's attention to aspects of a mathematical problem will likely be more productive if their initial focus is based on visual aspects or representations of the problem, and a reason is that their attentional processes are then mainly visual and will potentially more readily activate those brain regions now known to be involved in mathematical thought.

To be clear, the human brain is enormously complex, and we are at the very beginning of understanding how the brain processes, creates, discovers, and generally engages with mathematics. It is important therefore not to oversimplify the complexity of the issues involved. Presently we have hints and suggestions as to what, from an educational perspective, might work in accordance with recent neurological findings. As researchers move forward with neurological experiments we hope and expect that mathematics educators' research-based insights will contribute to the direction and interpretation of these experiments.

Interactions in mathematics classrooms can be enormously complex, with individual student brains engaging in mathematical thought, and communication between students and teacher taking place orally, in written language, and through specific mathematical signs and diagrams. A major area of ignorance is how individuals process these complex interactions, involving both language and non-language areas of their brains. For mathematics educators, seeking to find more productive ways of engaging students and assisting teachers, this, it seems to us, is a major point of connection between mathematics education and experimental neuroscience. As mathematics educators we need to understand the neuroscientific findings so we can relate them to our own experiences and, in turn, put sensible, productive questions to the neuroscientists.

**References**

Amalric, M., & Dehaene, S. (2016). Origins of the brain networks for advanced mathematics in expert mathematicians. *Proceedings of the National Academy of Sciences,* 113(18), 4909-4917.

Amalric, M., & Dehaene, S. (2018). Cortical circuits for mathematical knowledge: evidence for a major subdivision within the brain's semantic networks. *Philosophical Transactions of the Royal Society B: Biological Sciences*, 373(1740), 20160515.

Amalric, M., & Dehaene, S. (2019). A distinct cortical network for mathematical knowledge in the human brain. *NeuroImage*, 189, 19-31.

Amalric, M., Denghien, I., & Dehaene, S. (2018). On the role of visual experience in mathematical development: Evidence from blind mathematicians. *Developmental cognitive neuroscience*, *30*, 314-323.

Barlow, H. (2001). Redundancy reduction revisited. *Network: computation in neural systems*, *12*(3), 241-253.

Bettcher, B. M. (2017). Normal aging of the frontal lobes. In Miller, B. L., & Cummings, J. L. (Eds.)., *The Human Frontal Lobes. Functions and Disorders*, pp.343-356.

Bickmore-Brand, J. (Ed.). (1990). *Language in mathematics*. Australian Reading Association.

Bisley, J. W., & Goldberg, M. E. (2006). Neural correlates of attention and distractibility in the lateral intraparietal area. *Journal of neurophysiology*, *95*(3), 1696-1717.

Butterworth, B. (1999). *The Mathematical Brain*. Macmillan.

Camos, V., & Barrouillet, P. (2014). Attentional and non-attentional systems in the maintenance of verbal information in working memory: the executive and phonological loops. *Frontiers in human neuroscience*, *8*, 900.

Chapin, S. H., O'Connor, C., O'Connor, M. C., & Anderson, N. C. (2009). *Classroom discussions: Using math talk to help students learn, Grades K-6*. Math Solutions.
Clarke, D. J., Waywood, A., & Stephens, M. (1993). Probing the structure of mathematical writing. *Educational studies in mathematics*, 25(3), 235-250.
Cohen, R. A. (2014) *The Neuropsychology of Attention*. 2nd Edition. Springer.
Colby, C. L., Duhamel, J. R., & Goldberg, M. E. (1993). Ventral intraparietal area of the macaque: anatomic location and visual response properties. *Journal of neurophysiology*, *69*(3), 902-914.
Colby, C. L., & Goldberg, M. E. (1999). Space and attention in parietal cortex. *Annual review of neuroscience*, *22*(1), 319-349.
Culham, J. C., & Kanwisher, N. G. (2001). Neuroimaging of cognitive functions in human parietal cortex. *Current opinion in neurobiology*, *11*(2), 157-163.
Dehaene, S. (2009). Origins of mathematical intuitions: The case of arithmetic. *Annals of the New York Academy of Sciences*, 1156(1), 232-259.
Dehaene, S. (2010). The calculating brain. In Sousa, D. A. (Ed.), *Mind, brain, & education: Neuroscience implications for the classroom*, 179-198.
Dehaene, S. (2011). *The number sense: How the mind creates mathematics*. OUP USA.
Dehaene, S., Piazza, M., Pinel, P., & Cohen, L. (2003). Three parietal circuits for number processing. *Cognitive neuropsychology*, 20(3-6), 487-506.
Dehaene, S., Spelke, E., Pinel, P., Stanescu, R., & Tsivkin, S. (1999). Sources of mathematical thinking: Behavioral and brain-imaging evidence. *Science*, 284 (5416), 970-974.
DeJarnette, A. F., & González, G. (2016). Thematic analysis of students' talk while solving a real-world problem in geometry. *Linguistics and Education*, *35*, 37-49.
De Smedt, B., Ansari, D., Grabner, R. H., Hannula, M. M., Schneider, M., & Verschaffel, L. (2010) Cognitive neuroscience meets mathematics education. Educational Research Review, 5(1), 97-105.
De Smedt, B. & Grabner, R.H. (2015). Applications of neuroscience to mathematics education. In Kadosh, R.C. & Dowker, A. (Eds.), *Oxford library of psychology. The Oxford handbook of numerical cognition,* pp. 612–632. Oxford University Press.
De Smedt, B., & Verschaffel, L. (2010). Travelling down the road from cognitive neuroscience to education ... and back. *ZDM - The International Journal on Mathematics Education*, 42, 49-65.
Devlin, K. (2010). The mathematical brain. In Sousa, D. A. (Ed.), *Mind, brain, and education: Neuroscience Implications for the Classroom*, pp. 163-178.
DiCarlo, J. J., Zoccolan, D., & Rust, N. C. (2012). How does the brain solve visual object recognition?. *Neuron*, *73*(3), 415-434.
Dormal, V., & Pesenti, M. (2009). Common and specific contributions of the intraparietal sulci to numerosity and length processing. *Human brain mapping*, *30*(8), 2466-2476.
Ellerton, N. F., & Clarkson, P. C. (1996). Language factors in mathematics teaching and learning. In *International handbook of mathematics education* (pp. 987-1033). Springer, Dordrecht.
Gelman, R., & Butterworth, B. (2005). Number and language: how are they related?. *Trends in cognitive sciences*, *9* (1), 6-10.
Gottlieb, J. (2007). From thought to action: the parietal cortex as a bridge between perception, action, and cognition. *Neuron*, *53* (1), 9-16.
Grabner, R.H., Ansari, D., Schneider, M., De Smedt, B., Hannula, M.M. & Stern, E. (2010). Cognitive neuroscience and mathematics learning. *ZDM*, 42, 6.
Gray, E., & Tall, D. (2007). Abstraction as a natural process of mental compression. *Mathematics Education Research Journal*, *19* (2), 23-40.

Grotheer, M., Jeska, B., & Grill-Spector, K. (2018). A preference for mathematical processing outweighs the selectivity for Arabic numbers in the inferior temporal gyrus. *Neuroimage*, 175, 188-200.

Hawes, Z., Sokolowski, H. M., Ononye, C. B., & Ansari, D. (2019). Neural underpinnings of numerical and spatial cognition: An fMRI meta-analysis of brain regions associated with symbolic number, arithmetic, and mental rotation. *Neuroscience & Biobehavioral Reviews*, 103, 316-336.

Hunting, R.P. (Ed.) (1998) *Language Issues in Learning and Teaching Mathematics*. La Trobe University, Melbourne, australia

Hubbard, E. M., Piazza, M., Pinel, P., & Dehaene, S. (2005). Interactions between number and space in parietal cortex. *Nature Reviews Neuroscience*, *6*(6), 435-448.

Izard, V., & Dehaene, S. (2008). Calibrating the mental number line. *Cognition*, 106(3), 1221-1247.

Kadosh, R. C., & Walsh, V. (2009). Numerical representation in the parietal lobes: Abstract or not abstract?. *Behav Brain Sci*, 32(3-4), 313-328.

Kaplan, J. T., Gimbel, S. I., & Harris, S. (2016). Neural correlates of maintaining one's political beliefs in the face of counterevidence. *Scientific reports*, *6*, 39589.

Kaufmann, L., Kucian, K., & von Aster, M. (2015). Development of the numerical brain. In Kadosh, R.C. & Dowker, A. (Eds.), *Oxford library of psychology. The Oxford handbook of numerical cognition,* pp. 485–501. Oxford University Press.

Kieran, C. (2017). Cognitive neuroscience and algebra: Challenging some traditional beliefs. In *And the Rest is Just Algebra* (pp. 157-172). Springer, Cham.

Krueger, F., Spampinato, M. V., Pardini, M., Pajevic, S., Wood, J. N., Weiss, G. H., Landgraf, S. & Grafman, J. (2008). Integral calculus problem solving: an fMRI investigation. *Neuroreport*, 19(11), 1095.

LaBar, K. S., Gitelman, D. R., Parrish, T. B., & Mesulam, M. M. (1999). Neuroanatomic overlap of working memory and spatial attention networks: a functional MRI comparison within subjects. *Neuroimage*, *10*(6), 695-704.

Lafer-Sousa, R., & Conway, B. R. (2013). Parallel, multi-stage processing of colors, faces and shapes in macaque inferior temporal cortex. *Nature neuroscience*, *16*(12), 1870.

Lee, K., Lim, Z. Y., Yeong, S. H. M., Ng, S. F., Venkatraman, V., & Chee, M. W. L. (2007). Strategic differences in algebraic problem solving: Neuroanatomical correlates. *Brain Research*, 1155, 163–171.

McGowen, M. A., & Davis, G. E. (2019). Spectral analysis of concept maps of high and low gain undergraduate mathematics students. *The Journal of Mathematical Behavior*, 55, 100686.

Maruyama, M., Pallier, C., Jobert, A., Sigman, M. & Dehaen,e S. (2012). The cortical representation of simple mathematical expressions. *Neuroimage* 61, 1444 – 1460.

Menon, V. (2010). Developmental cognitive neuroscience of arithmetic: implications for learning and education. *ZDM*, *42*(6), 515-525.

Menon, V. (2015). Arithmetic in the child and adult brain. In Cohen Kadosh, R. & Dowker, A. (Eds.)*The Oxford Handbook of Numerical Cognition*, pp. 502-530. Oxford University Press.

Miyashita, Y. (1993). Inferior temporal cortex: where visual perception meets memory. *Annual review of neuroscience*, *16*(1), 245-263.

Monti, M. M., Parsons, L. M. & Osherson, D.N. (2012). Thought beyond language: neural dissociation of algebra and natural language. *Psychol. Sci*. 23, 914 – 922.

Morgan, C. (2006). What does social semiotics have to offer mathematics education research?. *Educational studies in mathematics*, *61*(1-2), 219-245.

Moschkovich, J. N. (2015). Academic literacy in mathematics for English learners. *The Journal of Mathematical Behavior*, *40*, 43-62.

O'Halloran, K. L. (2015). The language of learning mathematics: A multimodal perspective. *The Journal of Mathematical Behavior*, 40, 63-74.

Olshausen, B. A., & Field, D. J. (2000). Vision and the coding of natural images: The human brain may hold the secrets to the best image-compression algorithms. *American Scientist*, *88*(3), 238-245.

Peters, L., & De Smedt, B. (2018). Arithmetic in the developing brain: A review of brain imaging studies. *Developmental Cognitive Neuroscience*, *30*, 265-279.

Pinheiro-Chagas, P., Daitch, A., Parvizi, J., & Dehaene, S. (2018). Brain mechanisms of arithmetic: a crucial role for ventral temporal cortex. *Journal of cognitive neuroscience*, 30(12), 1757-1772.

Powell, S. R., Driver, M. K., Roberts, G., & Fall, A. M. (2017). An analysis of the mathematics vocabulary knowledge of third-and fifth-grade students: Connections to general vocabulary and mathematics computation. *Learning and Individual Differences*, *57*, 22-32.

Puglisi-Allegra, S., & Ventura, R. (2012). Prefrontal/accumbal catecholamine system processes high motivational salience. *Frontiers in behavioral neuroscience*, 6, 31.

Rickard, T. C., Romero, S. G., Basso, G., Wharton, C., Flitman, S., & Grafman, J. (2000). The calculating brain: an fMRI study. *Neuropsychologia*, 38(3), 325-335.

Sablé-Meyer, M., Caparos, S., van Kerkoerle, T., Amalric, M., & Dehaene, S. (2020). A signature of human uniqueness in the perception of geometric shapes. *PsyArxiv preprints*.

Schleppegrell, M. J. (2010). Language in mathematics teaching and learning: A research review. *Language and mathematics education: Multiple perspectives and directions for research*, 73-112.

Sellars, M. (2018). The Mathematical Brain. In *Numeracy in Authentic Contexts* (pp. 39-56). Springer, Singapore.

Shield, M., & Galbraith, P. (1998). The analysis of student expository writing in mathematics. *Educational Studies in Mathematics*, *36* (1), 29-52.

Spear, L. P. (2000). The adolescent brain and age-related behavioral manifestations. *Neuroscience & biobehavioral reviews*, *24* (4), 417-463.

Shum, J., Hermes, D., Foster, B. L., Dastjerdi, M., Rangarajan, V., Winawer, J., Miller, K. J. & Parvizi, J. (2013). A brain area for visual numerals. *Journal of Neuroscience*, *33* (16), 6709-6715.

Tall, D. (1978) Mathematical thinking and the brain. In *Proceedings of PME 2, Osnabrücker Schriften zür Mathematik*, pp. 333– 343.

Tall, D. (1994). The psychology of advanced mathematical thinking: Biological brain and mathematical mind. *In Proceedings of the 18th Conference of the International Group for the Psychology of Mathematics Education*, Vol. 1, pp. 33-39.

Tall, D. (2000). Biological brain, mathematical mind and computational computers. In *Proceedings of the 5th Asian Technology Conference in Mathematics* (pp. 3-20).

Tall, D. (2004). Thinking through three worlds of mathematics. In Høines, M. J. & Fuglestad, A. B. (Eds.) *Proceedings of the 28th Conference of the International Group for the Psychology of Mathematics Education*, Vol 4, pp. 281–288.

Tall, D. O. (2007). Developing a theory of mathematical growth. *ZDM*, 39 (1-2), 145-154.

Tall, D. O. (2019a). From biological brain to mathematical mind: The long-term evolution of mathematical thinking. In Danesi, M. (Ed.): *Interdisciplinary Perspectives on Math Cognition*, pp.1–28. Springer.

Tall, D. O. (2019b). Complementing supportive and problematic aspects of mathematics to resolve transgressions in long-term sense making. In *The Fourth Interdisciplinary Scientific Conference on Mathematical Transgressions*.

Tall, D. O. (2020a). Building long-term meaning in mathematical thinking: Aha! and Uh-Huh. https://bit.ly/370cmAT

Tall, D. O. (2020b). Making Sense of Mathematical Thinking over the Long Term: The Framework of Three Worlds of Mathematics and New Developments. *MINTUS: Beiträge zur mathematischen, naturwissenschaftlichen und technischen Bildung. Wiesbaden: Springer*.

Thompson, C. (2013). *Smarter than you think: How technology is changing our minds for the better*. Penguin

Thurston, W. P, (1990). Mathematical education. *Notices of the AMS*, 37, 844– 850.

Todd, J. J., & Marois, R. (2004). Capacity limit of visual short-term memory in human posterior parietal cortex. *Nature*, 428(6984), 751-754.

Vandierendonck, A. (2014). Symbiosis of executive and selective attention in working memory. *Frontiers in human neuroscience*, *8*, 588.

Verschaffel, L., Lehtinen, E., & Van Dooren, W. (2016). Neuroscientific studies of mathematical thinking and learning: a critical look from a mathematics education viewpoint. *ZDM*, *48*(3), 385-391.

Weber, K. (2001). Student difficulty in constructing proofs: The need for strategic knowledge. *Educational studies in mathematics*, *48*(1), 101-119.

Wilkinson, L. C. (2018). Teaching the language of mathematics: What the research tells us teachers need to know and do. *The Journal of Mathematical Behavior*, 51, 167-174.

Wilkinson, L. C. (2019). Learning language and mathematics: A perspective from Linguistics and Education. *Linguistics and Education*, *49*, 86-95.

**Declarations**

**Funding:** Not applicable

**Conflicts of interest/Competing interests:** Not applicable

**Availability of data and material:** Not applicable

**Code availability:** Not applicable

**Authors' contributions:** The authors contributed equally to the writing of the manuscript. All authors read and approved the final manuscript.

**Ethics approval** (include appropriate approvals or waivers) Not applicable

**Consent to participate** (include appropriate statements) Not applicable

**Consent for publication** (include appropriate statements) Not applicable